\documentclass[a4paper, 11pt]{article}

\usepackage{fullpage,amsmath,amssymb,amsthm,bm,url}
\usepackage[shortlabels]{enumitem}

\newcommand*{\prob}[1]{\mathbb P(#1)}
\newcommand*{\mean}[1]{\mathbb E(#1)}

\newcommand*{\floorfrac}[2]{\genfrac{\lfloor}{\rfloor}{}{}{#1}{#2}}
\newcommand*{\ceilfrac}[2]{\genfrac{\lceil}{\rceil}{}{}{#1}{#2}}

\newcommand{\ie}{i.e.\ }
\newcommand{\on}{\operatorname}
\newcommand{\ft}[1][t]{\mathcal F_{#1}}

\newtheorem{theorem}{Theorem} 
\newtheorem{lemma}[theorem]{Lemma}
\newtheorem{claim}[theorem]{Claim}

\theoremstyle{remark}

\title{Balanced two-type annihilation: mean-field asymptotics} 
	
\author{John Haslegrave\thanks{Department of Mathematics and Statistics, Lancaster University, UK. Supported by ERC Advanced Grant 883810.} \and 
	Peter Keevash\thanks{
		Mathematical Institute, University of Oxford, UK. Supported by ERC Advanced Grant 883810.}
}

\begin{document}\maketitle

\begin{abstract}We consider an interacting particle system where equal-sized populations  
	of two types of particles move by random walk steps on a graph,
	the two types  may have different speeds,
	and meetings of opposite-type particles result in annihilation. 
	The key quantity of interest is the expected extinction time. 
	Even for the mean-field setting of complete graphs, 
	the correct order of magnitude was not previously known. 
	Under essentially optimal assumptions on the starting configuration,
	we determine not only the order of magnitude but also the asymptotics:
	the expected extinction time on $K_{2n}$ is $(2+o(1))n\log n$,
	independently of the relative speeds of the two types.
	
	\medskip\noindent{\footnotesize\textbf{Keywords:} interacting particle system; reaction model; annihilation; random walk; complete graph\\
	\textbf{MSC2020 subject classifications:} 05C81; 82C22; 60K35.}
	\end{abstract}

\section{Introduction}

Interacting particle systems have been intensively studied since the 1970s
to model a variety of phenomena in Statistical Physics and Mathematical Biology,
such as spin systems, chemical reactions, population dynamics and the spread of infections.
An important example is the two-type annihilation model
for chemical reactions, studied classically on integer lattices $\mathbb{Z}^d$, 
where a series of papers by Bramson and Lebowitz \cite{BL90,BL91,BL01}
gives order of magnitude estimates for the site occupancy probability.
In this model, there are two types of particles and interactions only occur between
two particles of different type, which makes the analysis particularly difficult.
A single type model introduced by Erd\H{o}s and Ney \cite{EN74} is easier to analyse,
although even here the basic question of whether the origin would be occupied infinitely often was not solved in the original paper,
but was subsequently resolved in one dimension by Lootgieter \cite{Loo77} and in higher dimensions by Arratia \cite{Arr81,Arr83}. 
Besides these diffusive models, ballistic annihilation models have also been studied in the physics literature 
since the 1990s \cite{EF85,KRL95}, and more recently analysed rigorously in the mathematics literature \cite{HST21}.

Considering the real-world motivations for interacting particle systems, it is natural to consider systems
on finite graphs and ask for the asymptotics of key parameters when the number $n$ of vertices is large.
While a wide range of interacting particle systems have been well studied in a lattice setting (see e.g.\ \cite{Lig85} for an overview), the literature covering interacting particle systems on general, or random, finite graphs is less developed, and has mainly covered opinion-forming or information-spreading processes (see \cite{CdH24} for an overview), with one of the main models considered in this setting being the voter model, considered on general graphs by Donnelly and Welsh in the 1980s \cite{DW83}. Consequently there has been significant interest in the dual process of coalescing random walks on finite graphs (e.g.\ \cite{CEOR}).
Coalescence is simpler to analyse than annihilation due to methods based on duality and monotonicity.
A more general framework, with a variety of models including annihilation, was analysed by Cooper, Frieze and Radzik \cite{CFR09},
assuming that (a) the initial populations are polynomially small and randomly distributed,
(b) the underlying geometry is a random regular graph, and (c) the types move simultaneously with equal speeds.
By contrast, we will consider an annihilation model with individual rather than simultaneous movements
starting from an arbitrary spanning configuration, so we cannot avoid considering the so-called `big bang regime'
in which there is rapid action at the beginning of the process that is harder to analyse.

We refer to our companion paper \cite{expander} 
for a more detailed study of particle systems on expanders, 
with more history and new results for a variety of models, 
including the balanced two-type annihilation to be considered here;
our main result here will be a more precise asymptotic result in the mean-field setting.
Our precise formulation of the model is taken from a recent paper by Cristali, Jiang, Junge, Kassem, Sivakoff and York \cite{2-type}.
Fix a graph with $2n$ vertices, and initialise with one particle at each vertex, with $n$ particles being blue and $n$ red. 
The particles represent molecules of two reactants; when a red and a blue particle occupy the same vertex,
they react and are removed from the system, \ie they mutually annihilate. At each time step, 
with probability $p$ a red particle is chosen, uniformly at random from the remaining red particles,
or otherwise a uniformly random blue particle is chosen. Without loss of generality, we assume $p\leq 1/2$.
The chosen particle performs a simple random walk step. 
If it reaches a vertex with one or more particles of the opposite colour, 
it mutually annihilates with one such particle. This process almost surely 
eventually terminates with no particles remaining.

The key quantity of interest is the time taken for this to happen,
which we call the \emph{extinction time}. In analysing this we use standard asymptotic notation, with all limits applying as $n\to\infty$ unless otherwise stated: $f(n)=o(g(n))$ means $f(n)/g(n)\to 0$; $f(n)=\omega(g(n))$ means $f(n)/|g(n)|\to+\infty$; $f(n)=O(g(n))$ means $\limsup|f(n)/g(n)|<\infty$; and $f(n)=\Theta(g(n))$ means $0<\liminf |f(n)/g(n)|\leq\limsup |f(n)/g(n)|<\infty$.

The above model was considered by Cristali et al.~\cite{2-type} for the complete graph $K_{2n}$ (the mean-field case) and the star.
For these graphs, the one-type annihilation model is easy to analyse, but the two-type process is much harder, 
owing to the possibility of multiple particles of the same type occupying a single site. 
While they were able to give precise results on the extinction time for stars, for the mean-field case,
also assuming that $p$ is bounded away from zero uniformly in $n$,
they gave a lower bound of $2n\log n$ and an upper bound of $20n\log^2n/\log\log n$, 
thus leaving the open problem of showing that it is $\Theta(n\log n)$.
In \cite{expander}  we solved this problem for any $p$ 
on all regular graphs with sufficiently strong spectral expansion.
Here we use different methods to give the following more precise result in the mean-field setting.

\begin{theorem}\label{thm:complete}
	Let $p=p(n)\in[0,1/2]$. Consider an initial possibly random configuration on $K_{2n}$ having $n$ particles of each type, 
	with each site having at most one type of particle, where $n-A$ sites are initially occupied by red particles.
	Let $T$ be the extinction time. If $\mean{A}=o(pn\log n)$ then $\mean{T}  =(2+o(1))n\log n$ as $n\to\infty$.
\end{theorem}

Note that if $p=\omega(1/\log n)$, \ie the reds are not that slow, 
then Theorem \ref{thm:complete} applies to any starting configuration.
The assumption on $A$ cannot be significantly relaxed,
as the result is not true for all starting configurations when $p$ is of order $1/\log n$. 
To see this, notice that if all reds start on the same vertex $v$ 
then extinction requires at least $n$ occasions 
where either a red moves or a blue moves onto $v$. 
The number of these events in time $3n\log n$
is with high probability at most $3p n\log n+O(\log n)$,
so if $p=1/(4\log n)$ then extinction will not occur.

The idea of the proof is to control the distribution of the particles by martingale arguments, 
categorising it as good, intermediate or bad for each colour separately. 
If both distributions remain good or intermediate until very few particles remain, 
this will imply sufficiently good bounds on annihilation probabilities to give the desired result. 
We show that it is unlikely that random walk steps alone 
can cause a colour to cross the intermediate regime until late in the process. 
For blue particles, this also applies when we take annihilations into account, 
since they occur less frequently than blue walk steps, 
but for rarely-moving red particles it is possible that annihilations have a dominant effect. 
We thus bound both the time taken for the red distribution 
to recover from the effects of each annihilation, 
and the total number of annihilations that occur while the red distribution is bad. 
This part of the proof could be avoided if we assumed bounded speed ratios as in~\cite{2-type},
but this would be an unnecessary and moreover a strange assumption, 
as the limiting case when one colour is stationary seems intuitively easier!
The potential difficulty that must be avoided when $p$ is small but nonzero is that
the process may be slow to recover from a situation where due to blue annihilations 
the red particles have somehow become clustered on a few vertices.

\section{Proof}

Consider an initial possibly random configuration on $K_{2n}$ having $n$ particles of each type, 
with each site having at most one type of particle, where $n-A$ sites are initially occupied by red particles.
For simplicity we will add a loop at every vertex, so that each random walk step is equally likely to go to any vertex, 
including the current vertex. This produces a modified process that is just a time change by a factor $1+O(1/n)$,
so this makes no difference to our result. We also assume throughout that $n$ is sufficiently large.

At time $t$, we write $M_t$ for the number of particles of each colour,
$R_t$ for the number of sites with red particles, $B_t$ for the number of sites with blue particles, and
$\ft$ for the $\sigma$-algebra generated by all previous movements of particles.
When tracking $R_t$, we consider the effects of red and blue movements separately. To that end, we introduce the notation $R^*_t$, for $t\geq 1$, defined as follows: if a red particle moved at time $t$, then $R^*_t=R_t$; otherwise $R^*_t=R_{t-1}$. We can think of this as first choosing whether a red particle moves, if so moving it and resolving any annihilation that occurs, then updating the number of red-occupied sites to $R^*_t$, then moving a blue particle if appropriate, and finally updating the number of red-occupied sites to $R_t$. Since a blue movement can only affect the red particles via a single annihilation, we necessarily have $R_{t}\in \{R^*_{t}-1,R^*_{t}\}$. Note that if $R_t<R^*_t$ then a blue particle moved onto a lone red particle at time $t$.

We will show that with high probability $R_t$ and $B_t$ remain close to $M_t$ 
until almost all particles have been destroyed, except for short excursions. 
To quantify the closeness, we introduce a function $f(m)=\max(6m/n,1/w(n))$, 
where $w(n)=\omega(1)$, $w(n)^2=o(\log n)$ and $\mean{A}w(n)^2=o(pn\log n)$,
noting that such a function $w$ exists by our assumption on $\mean{A}$.
We say that blue is \emph{good} if $B_t(1+f(M_t))\geq M_t$, \emph{bad} if $B_t(1+2f(M_t))<M_t$, 
or  \emph{intermediate} otherwise, 
and likewise for red with respect to $R_t$. 

We will show that blue quickly becomes good, 
and thereafter is unlikely to become bad within the timeframe of interest.
On the other hand, the slower red particles are harder to control due to fluctuations
in their distribution caused by annihilations with the faster blue particles.
To analyse the red distribution, we introduce the red `level', $\ell_t$, which we define below.
We also maintain a stack of $\ell_t$ `thresholds', denoted $S_1,\ldots,S_{\ell_t}$ if $\ell_t>0$. The stack can only change in one of three ways: a new threshold is added when the level increases; the most recently added threshold is removed when the level decreases; or the entire stack is emptied when the level is reset to $0$. The initialisation and progression of the level and stack is as follows.
\begin{enumerate}[(a),nosep]
\item\label{level-1} If the initial configuration of red is good then $\ell_0=0$. Otherwise, $\ell_0=n-R_0$ and $S_i=n-i+1$ for $i=1,\ldots,\ell_0$.
\item\label{level-2} If red becomes good at time $t$ (\ie it is good at time $t$ but not time $t-1$), then we set $\ell_t=0$ (and empty the stack).
\item\label{level-3} If $R_t<R^*_t$ and red is not good, then we set $\ell_t=\ell_{t-1}+1$ and add $R^*_t$ to the stack.
\item\label{level-4} If $R_t=S_{\ell_{t-1}}$ then we set $\ell_t=\ell_{t-1}-1$ and remove $S_{\ell_{t-1}}$ from the stack.
\item\label{level-5} If $t\geq 1$ but none of the previous three cases apply, then $\ell_t=\ell_{t-1}$.
\end{enumerate}
We observe the following simple properties.
\begin{lemma}\label{lem:levels}For each time $t$, provided $\ell_t>0$ we have $S_1>\cdots>S_{\ell_t}>R_t$. Consequently, the next time $t'$ with $\ell_{t'}<\ell_t$ is also the next time such that either $R_{t'}=S_{\ell_t}$ or red is good at time $t'$. Furthermore, if red was good at time $t-1$, then red must be good at time $t'$.
\end{lemma}
\begin{proof}
We prove the first statement by induction; by \ref{level-1}, it holds for $t=0$. Suppose it holds at time $t-1$. If \ref{level-2} or \ref{level-4} apply, it clearly holds at time $t$. If \ref{level-3} applies, then $R_t<R^*_{t}=S_{\ell_t}$ and blue must have moved at time $t$, meaning $R^*_t=R_{t-1}$. Thus, provided $\ell_t>1$, we also have $S_{\ell_t}=R_{t-1}<S_{\ell_t-1}$. Finally, if \ref{level-5} applies then $R_t\leq R_{t-1}+1\leq S_{\ell_{t-1}}$, and since \ref{level-4} does not apply at least one inequality must be strict.

Note that the next time $t'$ with $\ell_{t'}<\ell_t$ coincides with the time that $S_{\ell_t}$ is removed from the stack. This can only happen if either red is good at time $t'$ or $R_{t'}=S_{\ell_t}$. If the former holds then $\ell_{t'}=0$ by \ref{level-2}. If the latter holds then $S_{\ell_t}$ can no longer be on the stack, since this would contradict the first statement. Thus the second statement follows.

To prove the third statement, note that red being good at time $t-1$ implies $\ell_{t-1}=0$, and so we must have $\ell_t=1$, $R^*_t=R_{t-1}$ and $R_t=R^*_t-1$, giving $S_{\ell_t}=R_{t-1}$. It is easy to check that the function $x\mapsto \frac{x}{1+f(x)}$ is increasing, and so $R_{t'}=S_{\ell_t}$ would imply red is good at $t'$.
\end{proof}

We define a stopping time $\tau$ as the first time that one of the following occurs:
\begin{enumerate}[(i),nosep]
	\item $M_t = \lfloor\log^2 n\rfloor$;\label{success}
	\item blue is bad but has previously been good;\label{f1}
	\item red is bad at level $0$;\label{f2} or
	\item $t = 3n^2$. \label{f3}
\end{enumerate}

Our estimate for the extinction time $T$ will be
$T \leq T_{\on{blue}} + T_{\on{red}} + T_{\on{late}} + T_{\on{ok}}$,
where
\begin{itemize}[nosep]
	\item $T_{\on{blue}}$ is the number of steps before blue first becomes good,
	\item $T_{\on{red}}$ is the number of steps when red is above level $0$, 
	\item $T_{\on{late}}$ is the number of steps after $\tau$ with red at level $0$, and
	\item $T_{\on{ok}}$  is the number of steps before $\tau$ where neither colour is bad.
\end{itemize}
To see that this estimate is valid, we first write $T_{\on{red}}=T_{\on{red}}^{(1)}+T_{\on{red}}^{(2)}$, where the term $T_{\on{red}}^{(1)}$ counts the number of steps before $\tau$ when red is above level $0$, and $T_{\on{red}}^{(2)}$ counts the number of steps after $\tau$ when red is above level $0$. Now consider any time step before $\tau$. If blue has not yet become good, this time step is counted by $T_{\on{blue}}$. If red is above level $0$, this time step is counted by $T_{\on{red}}^{(1)}$. If neither of these apply then blue has previously become good and red is at level $0$. By definition of $\tau$, (ii) and (iii) have not yet occurred, meaning that blue is still good and red is not bad, so this time step is counted by $T_{\on{ok}}$. Thus 
\begin{equation}\label{t-bound-1}\tau\leq T_{\on{blue}} + T_{\on{red}}^{(1)} + T_{\on{ok}}.\end{equation}
Finally, consider any time step after $\tau$. This time step is counted by either $T_{\on{red}}^{(2)}$ or $T_{\on{late}}$, depending on the level. Therefore $T-\tau= T_{\on{red}}^{(2)}+T_{\on{late}}$, which together with \eqref{t-bound-1} gives the required bound.

We start our analysis with the following simple tail bound on the runtime.
\begin{lemma}\label{quick-bound}
	From any starting position with $m\leq n$ particles of each type, 
	$\mean{T}\leq 2mn$ and $\prob{T>3n^2}\leq e^{-\Theta(n)}$.
\end{lemma}
\begin{proof}
	At each step, conditional on the previous evolution, an annihilation occurs with probability at least $\frac{1}{2n}$,
	so the expected time between consecutive annihilations is at most $2n$. Summing gives the bound for $\mean{T}$.
	The probability of needing  more than $3n^2$ steps for $m\leq n$ annihilations 
	is at most $\prob{\on{Bin}(3n^2,\frac{1}{2n})<n}$, so the probability bound holds by Chernoff.
\end{proof}

Next we consider the main term in the estimate for $T$.
\begin{lemma}\label{Tok} We have $\mean{T_{\on{ok}}}\leq 2n\log n+100 w(n)^{-1} n\log n$. \end{lemma}
\begin{proof}
	Each step when neither colour is bad results in an annihilation with probability at least
	\[\frac{p{B_t}+(1-p){R_t}}{2n}\geq \frac{M_t}{2n(1+2f(M_t))}.\]
	Thus we bound the expected number of steps 
	with neither colour bad before the next annihilation
	by  $26n/M_t$ while $M_t\geq n/6w(n)$ (so $f(M_t)\leq 6$)
	or by $(1+2/w(n))2n/M_t$ while $M_t\leq n/6w(n)$ (so $f(M_t)=1/w(n)$).
	Summing and using $\sum_{m=k}^n1/m\leq1+\log n-\log k$ for all $n>k\geq 1$, 
	\begin{align*}\mean{T_{\on{ok}}}&\leq(1+2/w(n))\sum_{m=1}^{n}2n/m+\sum_{m=\lceil n/6w(n)\rceil}^n 24n/m\\
		&\leq 2n\log n+100w(n)^{-1} n\log n,\end{align*}
	since $w(n)\leq \sqrt{\log n}$ and $n$ is sufficiently large.
\end{proof}

The evolution of $B_t$ may be bounded as follows, writing $q:=1-p$. We have
\begin{equation}\prob{B_{t+1}=B_t+1\mid\ft}\geq q\frac{M_t-B_t}{M_t}\cdot\frac{2n-B_t-R_t}{2n},\label{bt-1}\end{equation}
since this occurs when a blue particle at a site occupied by at least one other blue particle, of which there are at least $M_t-B_t$, moves to an empty site.
Similarly, we have
\begin{equation}\prob{B_{t+1}=B_t-1\mid\ft}\leq q\frac{B_t}{M_t}\cdot\frac{B_t+R_t}{2n}+p\frac{B_t}{2n}\label{bt-2},\end{equation}
since for this to happen either the only blue particle at some site must move to a non-empty site, or a red particle must move to a site occupied by exactly one blue particle. In all other cases, $B_{t+1}=B_t$.

Recall that we distinguish changes in $R_t$ owing to red movements from those owing to blue movements,
which may be dominant if blue is much faster than red. Not including the effect of blue movements, the corresponding transitions are
\begin{equation}\prob{R^*_{t+1}=R_t+1\mid\ft}\geq p\frac{M_t-R_t}{M_t}\cdot\frac{2n-R_t-B_t}{2n}\label{rt-1}\end{equation}
and
\begin{equation}\prob{R^*_{t+1}=R_t-1\mid\ft}\leq p\frac{R_t}{M_t}\cdot\frac{R_t+B_t}{2n}.\label{rt-2}\end{equation}

We will see that when a colour is not good then increments are significantly more likely, giving a strong self-correcting effect.
To describe this, we need the following simple facts about biased random walks.
\begin{lemma}\label{bias-walk}
	Let $X_t$ be an integer-valued discrete-time stochastic process adapted to a filtration $(F_t)_{t\geq 0}$ and stopped at some stopping time $S$. Suppose, for all $t\leq S$, we have
	\begin{itemize}
		\item $X_{t+1}\in\{X_t-1,X_t,X_t+1\}$; and
		\item $\prob{X_{t+1}=X_t-1\mid F_t}\leq \prob{X_{t+1}=X_t+1\mid F_t}/2$.
	\end{itemize}
	Then, given $X_0=k-1$, and conditioned on $F_0$, the probability of reaching $0$ before $k$ is at most $2^{1-k}$. If additionally there is some $0<\alpha\leq 2/3$ such that $\prob{X_{t+1}=X_t+1\mid F_t}\geq \alpha$ for all $t$, then the expected time to reach $k$ (or stop without having hit $k$) is at most $2\alpha^{-1}$.
\end{lemma}
\begin{proof}
	For the first statement, it is sufficient to prove for every finite time $t'$ that the probability of reaching $0$ before reaching $k$ and before time $t'$ is at most $2^{1-k}$. Note that $2^{k-X_t}$ is a supermartingale, and the first time $T$ for which $X_T\in\{0,k\}$ or $T=S\wedge t'$ is a stopping time. 
	Since $2^{k-X_{t\wedge T}}$ is bounded, the Optional Stopping Theorem
	gives $\mean{2^{k-X_T}\mid F_0}\leq 2^{k-X_0}=2$. Since $2^{k-X_T}\geq 1$ we obtain $\prob{X_T>0\mid F_0}+2^k\prob{X_T=0\mid F_0}\leq2$, 
	and therefore $\prob{X_T=0\mid F_0}\leq 2^{1-k}-2^{-k}\prob{X_T>0\mid F_0}<2^{1-k}$.
	
	For the second statement, let $Z^+_t$ be the number of $s<t$ such that $X_{s+1}=X_s+1$, and similarly let $Z^-_t$ be the number of $s<t$ such that $X_{s+1}=X_s-1$. Let $T'$ be the first time that $X_{T'}=k$ or $T'=S$. Note that $T'$ is a stopping time of finite expectation, and that $Z^+_t-2Z^-_t$ is a submartingale. Thus the Optional Stopping Theorem gives $\mean{Z^+_{T'}\mid F_0}\geq 2\mean{Z^-_{T'}\mid F_0}$, and by definition of $T'$ we have $Z^+_{T'}-Z^-_{T'}\leq 1$. Thus $\mean{Z^+_{T'}\mid F_0}-\mean{Z^+_{T'}\mid F_0}/2\leq 1$, so $\mean{Z^+_{T'}\mid F_0}\leq 2$. Finally, $Z^+_t-\alpha t$ is a submartingale, so $\mean{T'\mid F_0}\leq\alpha^{-1}\mean{Z^+_{T'}\mid F_0}\leq 2\alpha^{-1}$. 
\end{proof}

We next establish the necessary bias in our walks. 
\begin{lemma}\label{bias}
	Provided that the blue distribution is not good, we have 
	\begin{equation}\frac{\prob{B_{t+1}=B_t+1\mid\ft}}{\prob{B_{t+1}=B_t-1\mid\ft}}\geq 2,\label{b-bound}\end{equation}
	and similarly, provided that the red distribution is not good, we have
	\begin{equation}\frac{\prob{R^*_{t+1}=R_t+1\mid\ft}}{\prob{R^*_{t+1}=R_t-1\mid\ft}}\geq 2.\label{r-bound}\end{equation}
\end{lemma}
\begin{proof}We prove \eqref{b-bound}. Since the bound in \eqref{rt-1} matches that of \eqref{bt-1}, after exchanging $R_t$ and $B_t$, and $p$ and $q$, but that of \eqref{rt-2} is stronger than \eqref{bt-2} (after the same exchanges), the proof of \eqref{r-bound} only differs in having more room to spare.
	Since $R_t\leq M_t$ and $p\leq q$, we have, using \eqref{bt-1}, 
	\[\prob{B_{t+1}=B_t+1\mid\ft}\geq q\frac{M_t-B_t}{M_t}\cdot\frac{2n-B_t-M_t}{2n}\]
	and, using \eqref{bt-2},
	\[\prob{B_{t+1}=B_t-1\mid\ft}\leq q\frac{B_t}{M_t}\cdot\frac{B_t+M_t}{2n}+p\frac{B_t}{2n}\leq q\frac{B_t}{M_t}\cdot\frac{B_t+2M_t}{2n},\]
	and so it is sufficient to verify that
	\[(M_t-B_t)(2n-B_t-M_t)\geq 2B_t(B_t+2M_t),\]
	or equivalently that
	\begin{equation}\label{blue-bias}2n(M_t-B_t)\geq M_t^2+4M_tB_t+B_t^2.\end{equation}
	First suppose $M_t\geq 7B_t$. Then, since $n\geq M_t$, the LHS of \eqref{blue-bias} is at least $\frac{12}{7}M_t^2$ 
	and the RHS is at most $\frac{78}{49}M_t^2$, so \eqref{blue-bias} holds. 
	Now suppose $7B_t\geq M_t$. As blue is not good, we have $M_t \ge (1+6M_t/n)B_t$.
	Then the LHS of \eqref{blue-bias} is at least $12M_tB_t$ and the RHS is at most $11M_tB_t$, so again \eqref{blue-bias} holds.
\end{proof}

We deduce a suitable bound on the number of steps before blue first becomes good.
\begin{lemma}\label{Tblue} We have $\mean{T_{\on{blue}}}\leq 2n$.\end{lemma}
\begin{proof}We show that with probability $1-\exp(-\Theta(n))$, blue becomes good before $5n/9$ annihilations have occurred and within $\lceil 6nq^{-1}/7\rceil<2n$ time steps. The result then follows since we can bound the expected additional number of steps required when this fails by $2n^2\exp(-\Theta(n))=o(1)$, using Lemma \ref{quick-bound}.
	
Note that blue is good whenever $B_t\geq n/7$. Provided $B_t<n/7$ and $M_t>4n/9$, we can calculate using \eqref{bt-1} and $R_t\leq M_t$ that
\[\prob{B_{t+1}=B_t+1\mid\ft}\geq q\min_{x\in[4/9,1]}\left(\frac{x-1/7}{x}\cdot\frac {13/7-x}{2}\right)=18q/49.\]
Similarly, since $p\leq q$, we have from  \eqref{bt-2} that
\[\prob{B_{t+1}=B_t-1\mid\ft}\leq q\max_{x\in[4/9,1]}\left(\frac{1/7+x}{14x}+\frac {1}{14}\right)=65q/392.\]
Thirdly, the probability of an annihilation under the same assumptions is at most
\[q\frac {M_t}{2n}+p\frac{B_t}{2n}\leq 4q/7.\]
Now run the process until one of the following occurs: $M_t\leq 4n/9$, $B_t\geq n/7$ or $\lceil 6nq^{-1}/7\rceil$ time steps have occurred. Note that the number of annihilations up to this point is dominated by $\on{Bin}(\lceil 6nq^{-1}/7\rceil,4q/7)$, which by a standard Chernoff bound has exponentially small probability of exceeding $5n/9$. Similarly, the number of times $B_t$ decreases up to this point has exponentially small probability of exceeding $n/7$. Since $\on{Bin}(\lceil 6nq^{-1}/7\rceil,18q/49)<2n/7$ with exponentially small probability, $B_t$ reaches $n/7$ before time $\lceil 6nq^{-1}/7\rceil$ with suitably high probability.
\end{proof}

Recall that the stopping time $\tau$ is the first time when one of four possible conditions holds. The `success' event is that \ref{success} occurs at time $\tau$, \ie $M_\tau = \lfloor\log^2 n\rfloor$,
and the `failure' event, denoted $F$, is that one of \ref{f1}, \ref{f2} or \ref{f3} occurs before \ref{success}. We now bound the probability of failure.

\begin{lemma} \label{F} We have $\prob{F}\leq 3n^{-3}$.\end{lemma}
\begin{proof}
	Note that $\prob{F}\leq \prob{F_{\on{blue}}}+\prob{F_{\on{red}}}+\prob{F_{\on{time}}}$, where $F_{\on{blue}}$ is the event that \ref{f1} occurs (that is, blue is bad but has previously been good) at time $\tau$, and likewise $F_{\on{red}}$ that \ref{f2} does (that is, red is bad at level $0$) and $F_{\on{time}}$ that \ref{f3} does (that is, $\tau=3n^2$).
	
	Suppose that $F_{\on{blue}}$ or $F_{\on{red}}$ occurs. We say that the \textit{genesis} of such an event is the last time $\tau^*$ prior to $\tau$ that the relevant colour was good. For $F_{\on{blue}}$, the existence of such a time is immediate from \ref{f1}; for $F_{\on{red}}$ note that \ref{level-1} ensures that red will be good when the level first reaches $0$.
	
	We first bound, for each fixed time $t^*<3n^2$, the probability that blue is good at $t^*$, not good at $t^*+1$, and thereafter becomes bad before becoming good. Denote this event by $E_{t^*,\on{blue}}$, and let $E'_{t^*,\on{blue}}$ be the event that blue is good at $t^*$ but not good at $t^*+1$. Note that, if $F_{\on{blue}}$ occurs with genesis $\tau^*$, then $E_{\tau^*,\on{blue}}$ occurs (and $\tau^*<\tau\leq 3n^2$), and so taking a union bound across all $t^*<3n^2$ will bound $\prob{F_{\on{blue}}}$.
				
	Clearly $\prob{E_{t^*,\on{blue}}}\leq \prob{E_{t^*,\on{blue}}\mid E'_{t^*,\on{blue}}}$, so we bound the latter; note that $E'_{t^*,\on{blue}}$ is $\mathcal F_{t^*+1}$-measurable. Let $X_i=B_{t^*+1+i}$ for $i=0,\ldots,S$, where $S$ is the smallest value such that blue is either good or bad at $t^*+1+S$.  We clearly have $X_{i}\in\{X_{i}-1,X_{i},X_{i}+1\}$ for all $0\leq i<S$. Furthermore, since blue is not good for all $t^*+1\leq t<S$, \eqref{b-bound} of Lemma \ref{bias} applies. Thus both conditions of Lemma \ref{bias-walk} are met with $F_t=\mathcal F_{t+t^*+1}$.
	
	Note that $B_{t^*}\geq \frac{M_{t^*}}{1+f(M_{t^*})}$ but $X_0=B_{t^*+1}<\frac{M_{t^*+1}}{1+f(M_{t^*+1})}$, which together with $B_{t^*+1}\geq B_{t^*}-1$ and the fact that $\frac{M_t}{1+f(M_t)}$ is non-increasing implies $X_{0}=\ceilfrac{M_{t^*}}{1+f(M_{t^*})}-1$.  In order for $E_{t^*,\on{blue}}$ to occur, $X_i$ must reach $\floorfrac{M_{t^*}}{1+2f(M_{t^*})}$ before $\ceilfrac{M_{t^*}}{1+f(M_{t^*})}$. Lemma \ref{bias-walk} gives that this occurs with probability at most $2^{1-z}$, where $z=\ceilfrac{M_{t^*}}{1+f(M_{t^*})}-\floorfrac{M_{t^*}}{1+2f(M_{t^*})}$.
		
	Next we consider $F_{\on{red}}$. For each fixed time $t^*<3n^2$, denote by $E_{t^*,\on{red}}$ the event that red is good at $t^*$, red is not good at time $t^*+1$, and that red becomes bad at level $0$ before it becomes good, and let $E'_{t^*,\on{red}}$ be the event that red is good at $t^*$ but not at time $t^*+1$, and $\ell_{t^*+1}=0$. Note that if  $\ell_{t^*+1}>0$ then, by Lemma \ref{lem:levels}, red is good at the next time $t'$ that $\ell_{t'}=0$, and so $E_{t^*,\on{red}}$ does not occur. Again, if $F_{\on{red}}$ occurs with genesis $\tau^*$, then $\tau^*<3n^2$ and $E_{\tau^*,\on{red}}$ and $E'_{\tau^*,\on{red}}$ both occur, so $\prob{F_{\on{red}}}\leq\sum_{t^*<3n^2}\prob{E_{t^*,\on{red}}\mid E'_{t^*,\on{red}}}$. 
	
	In order to bound $\prob{E_{t^*,\on{red}}\mid E'_{t^*,\on{red}}}$, condition on the latter event and let $t^*+1=c_0<c_1<\cdots$ be the times after $t^*$ such that $\ell_{c_i}=0$, let $S$ be the smallest value such that red is either good or bad at $c_S$ and let $X_i=R_{c_i}$ for $i=0,1,\ldots, S$. In other words, we remove any excursions above level $0$, each of which returns $R_t$ to the value before the excursion. This ensures $X_{j+1}=R^*_{c_{j}+1}$ for each $0\leq j<S$. In particular, $X_{j+1}\in \{R_{c_j}-1, R_{c_j}, R_{c_j}+1\}$, and, since red is not good, \eqref{r-bound} of Lemma \ref{bias} applies. Thus again both conditions of Lemma \ref{bias-walk} are met with $F_t=\mathcal F_{c_t}$. By the same argument as before, we have $X_{0}=\ceilfrac{M_{t^*}}{1+f(M_{t^*})}-1$, and in order for $E_{t^*,\on{red}}$ to occur, $X_i$ must reach $\floorfrac{M_{t^*}}{1+2f(M_{t^*})}$ before $\ceilfrac{M_{t^*}}{1+f(M_{t^*})}$, which by Lemma \ref{bias-walk} has probability at most $2^{1-z}$.
	
	Note that $\frac{x}{(1+x)(1+2x)}$ is increasing for $0\leq x\leq 2^{-1/2}$ and decreasing for $x\geq 2^{-1/2}$. 
	Also, we have $1/w(n)\leq f(M_t)\leq 6$ and $M_{t^*}\geq \log^2 n$. Consequently, we have
	\begin{align*}\ceilfrac{M_{t^*}}{1+f(M_{t^*})}-\floorfrac{M_{t^*}}{1+2f(M_{t^*})}&\geq \frac{M_{t^*}f(M_{t^*})}{(1+f(M_{t^*}))(1+2f(M_{t^*}))}\\
		&\geq\min(M_{t^*}/(4w(n)),6M_{t^*}/91)\geq 10\log n.
	\end{align*}
	Thus $2^{1-z}\leq 2^{1-10\log n}\leq n^{-6}$. By a union bound over all $t^*<3n^2$ and choices of blue or red, we have $\prob{F_{\on{blue}}}+\prob{F_{\on{red}}}\leq 2n^{-3}$.
	Finally, the probability of failure due to taking too many steps 
	is bounded by $\prob{T\geq 3n^2} \le n^{-3}$ by Lemma \ref{quick-bound}.
\end{proof}

Next we control the time spent with red not at level zero.
\begin{lemma}\label{Tred} We have $\mean{T_{\on{red}}}\leq 12nw(n)^2+6p^{-1}\mean{A}w(n)^2$.\end{lemma}
\begin{proof}
	We break the calculation into the following two claims, 
	where the first bounds the contribution from each increase in level, 
	and the second the expected number of increases in level. Here we consider that the level increases to $i$ at time $t^*$ if either $t^*>0$, $\ell_{t^*-1}=i-1$ and $\ell_{t^*}=i$, or $t^*=0$ and $\ell_0\geq i\geq 1$. Thus we consider that the process starts with $\ell_0$ simultaneous increases.
	\begin{claim}\label{clm:up-down}
		Suppose the level increases to $i$ at time $t^*$. Let $t'$ be the next time with $\ell_{t'}<i$. Let $C=\{t:t^*\leq t\leq t'\text{ and }\ell_t=i\}$. Then $\mean{|C|\mid \ft[t^*]}\leq 5p^{-1}(1+w(n))$.
	\end{claim}
	\begin{proof}
		We show the slightly stronger bound that $\mean{|C|\mid \ft[t^*], C\neq\varnothing}\leq 5p^{-1}(1+w(n))$, which only differs if $t^*=0$, since otherwise $t^*\in C$.
		Write $C=\{c_1,\ldots,c_s\}$, where $c_1<c_2\cdots<c_s<t'$. 		
		As in the proof of Lemma \ref{F}, we consider the sequence $R_{c_1},\ldots,R_{c_s}$, \ie we remove any excursions above level $i$, each of which returns $R_t$ to the value before the excursion. This ensures $R_{c_j+1}^*=R_{c_{j+1}}$ for each $j<s$.
		
		Note that, using Lemma \ref{lem:levels}, the level will fall below $i$ if $R_t$ exceeds $R_{c_1}$: if $t^*>0$ then $c_1=t^*$ and this follows by \ref{level-3} and \ref{level-4}, whereas if $t^*=0$ then it follows from \ref{level-1} and \ref{level-4}. 
		As red is not good at $c_j$ for any $j\leq s$, we have
		\begin{align*} \prob{R_{c_j+1}^*=R_{c_j}+1\mid \ft[c_j]}&\geq p\frac{M_{c_j}-R_{c_j}}{M_{c_j}}\cdot\frac{2n-R_{c_j}-B_{c_j}}{2n}\\
		&\geq p\frac{f(M_{c_j})}{1+f(M_{c_j})}\cdot\frac37\geq \frac{3p}{7}\cdot\frac{1}{1+w(n)}.\end{align*}
		By Lemma \ref{bias}, the conditions of Lemma \ref{bias-walk} are satisfied with $\alpha=\frac{3p}{7}\cdot\frac{1}{1+w(n)}$ and thus the expected time for this sequence to exceed its initial value is at most $14p^{-1}(1+w(n))/3$. 
	\end{proof}
	\begin{claim}\label{bad-moves}
		The expected number of increases in level is at most $(2pn+\mean{A})w(n)$.
	\end{claim}
	\begin{proof}
		Recall that $A=n-R_0$ and $q=1-p$. Let $L$ be the number of remaining red particles when the level first increases.
		It suffices to show $\mean{L \mid A=a} \le (2pn+a)w(n)$ for each $a$ such that $\prob{A=a}>0$.
		This is easy to see if red is not initially good, as then $a \ge 6n/7$ and $L=n<aw(n)$.
		
		Now suppose red is initially good. We say that a `bad move' by red is a movement of a red particle 
		onto another site occupied by red. We write $\Xi_t$ for the number of bad moves by time $t$.
		As $M_t-R_t$ can only increase via bad moves, we have $M_t-R_t \le a + \Xi_t \le a + \Xi_T$.
		When the level increases red is not good, so $M_t - R_t \ge R_t f(R_t) \ge R_t/w(n)$.
		We deduce $L \le (a + \Xi_T)w(n)$. 
		
		It remains to bound $\Xi_T$. We note that a bad move at time $t$ has probability at most $p\frac{R_t}{2n}$, 
		whereas the probability of an annihilation is at least $q\frac{R_t}{2n}$, 
		since this is the probability of a blue particle moving onto a site occupied by red particles. 
		Thus $pC_t-q\Xi_t$ is a submartingale, where $C_t$ is the number of annihilations by time $t$. 
		At the stopping time $T$, since $C_T=n$, we have $\mean{\Xi_T}\leq np/q\leq 2pn$, as required.
	\end{proof}
	Now we deduce the lemma from the following calculation combining these two claims:
	\begin{align*}\mean{T_{\on{red}}}&\leq(2pn+\mean{A})w(n)\cdot 5p^{-1}(1+w(n))\\
		&=(10n+5\mean{A}p^{-1})(w(n)+w(n)^2)\\
		&\leq 12nw(n)^2+6p^{-1}\mean{A}w(n)^2. \qedhere\end{align*}
\end{proof}

It remains to bound the time for completion once only $\log^2 n$ particles remain.

\begin{lemma} \label{Tlate} We have $\mean{T_{\on{late}}}\leq 26n\log\log n$.\end{lemma}
\begin{proof}
	Recalling that $F$ is the failure event, we write
	$\mean{T_{\on{late}}} = \prob{F}\mean{T_{\on{late}} \mid F} + \mean{T_{\on{late}}  1_{\neg F}}$.
	By Lemma \ref{F}, we have $\prob{F}\leq 3n^{-3}$. By Lemma \ref{quick-bound}, we have $\mean{T_{\on{late}}\mid \mathcal F_{\tau}}\leq 2n^2$. Since $F$ is $\mathcal F_{\tau}$-measurable, also $\mean{T_{\on{late}} \mid F}\leq 2n^2$ and so the first term is at most $6/n$.
	For the second, we assume $F$ does not occur, and note that
	$T_{\on{late}}$ only counts steps $t$ after $\tau$ with red at level $0$, 
	so red must be good at such $t$. In particular, $R_t\geq M_t/2$, 
	so the probability of an annihilation is at least $p\frac{M_t}{4n}\geq \frac{M_t}{8n}$.
	Thus this contribution has expected size at most
	$\sum_{m=1}^{\lfloor\log^2 n\rfloor}8n/m\leq 16n\log\log n+8n$.
\end{proof}

Combining the bounds from Lemmas \ref{Tok}, \ref{Tblue}, \ref{Tred} and \ref{Tlate}, using $\mean{A}w(n)^2=o(pn\log n)$ gives
\begin{align*}
	\mean{T}& \le 2n\log n+100 w(n)^{-1} n\log n + 2n
	+  12nw(n)^2 + 6p^{-1}\mean{A}w(n)^2 + 26n\log\log n \\
	& = (2+o(1))n\log n, \quad \text{as required to prove Theorem \ref{thm:complete}.}
\end{align*}

As a variant form of Theorem \ref{thm:complete},
one can obtain an explicit estimate for the second order term
under the stronger assumption $\mean{A}\leq pn$.
Indeed, then applying the argument with $w(n)=\log^{1/3}n$ we obtain
$\mean{T} \le 2n\log n + O(n\log^{2/3} n)$.
This form of the second order term seems to be an artefact
of the proof; we do not think it is optimal.

\section{Concluding remarks}

In a companion paper \cite{expander} we consider a more general class of interacting particle systems on expanders,
for which we obtain similar $\Theta(n\log n)$ bounds on the equilibrium time. Our methods there are entirely different,
and while they apply in much greater generality, they only reveal the order of magnitude, not the asymptotics.
It may be interesting to extend the methods of this paper to obtain asymptotic results for other graphs
with simple structure, such as the complete bipartite graph $K_{n,n}$. This question conceals another, more basic question: 
does $K_{n,n}$ have a well-defined asymptotic extinction time from any starting distribution, 
or do some initial distributions reach extinction much faster than others?

For $K_{n,n}$ we believe that the starting distribution does not make much difference.
We give some evidence for this by sketching a proof in the case of stationary reds.
In this case the extinction time does not depend on the order in which blue particles
are chosen to move (the `abelian sandpile' property, see \cite{expander}).
Thus we may consider blues in turn, moving each one until it hits a red. 
Each time taken is a random variable of the following form: 
count the number of independent trials required for success, 
where odd trials have some success probability $p_1$ 
and even trials have some success probability $p_2$. 
While $p_1$ and $p_2$ depend on the arrangement of remaining reds, 
$p_1+p_2=m/n$ does not, where $m$ is the number of reds remaining. 
Now this is $2X-Y$ where $X$ is a geometric random variable with success probability $p_1+p_2-p_1p_2$, 
\ie the probability that one of any given pair of trials succeeds, 
and $Y$ is an independent Bernoulli random variable with probability $p_1/(p_1+p_2-p_1p_2)$. 
In particular, $\mean{2X-Y}$ is minimised when $p_2=0$, at $2n/m-1$. 
Assuming $m<n/\log n$ we have $p_1p_2<(p_1+p_2)^2 /4<(p_1+p_2)/(4 \log n)$,
so $\mean{2X-Y} \le (2+1/\log n)m/n$, and the sum of all terms with $m<n/\log n$ is at most $2n\log n-2n\log\log n+O(n)$. 
The sum of all other terms is increased by at most a constant factor, 
and so has sum $O(n\log\log n)$. Thus for stationary reds on $K_{n,n}$ 
we have $\mean{T}=(2+o(1))n\log n$ for any starting distribution. 

\section*{Acknowledgements}
	Both authors were supported by ERC Advanced Grant 883810. We are grateful to the anonymous referees for their help in improving the exposition.

\end{document}